\documentclass[12pt]{article}

\usepackage{amsmath, epsfig, cite}
\usepackage{amssymb}
\usepackage{amsfonts}
\usepackage{latexsym}
\usepackage{float}
\usepackage{color}

\newtheorem{thm}{Theorem}[section]

\newtheorem{conj}[thm]{Conjecture}

\numberwithin{equation}{section}

\begin{document}

\begin{center}
{\large\bf Two new $q$-supercongruences arising from Carlitz's identity}
\end{center}

\vskip 2mm \centerline{Ji-Cai Liu and Wei-Wei Qi}
\begin{center}
{\footnotesize Department of Mathematics, Wenzhou University, Wenzhou 325035, PR China\\[5pt]
{\tt jcliu2016@gmail.com, wwqi2022@foxmail.com} \\[10pt]
}
\end{center}

\vskip 0.7cm \noindent{\bf Abstract.}
From Carlitz's identity, we deduce two new $q$-supercongruences modulo the square of a cyclotomic polynomial,
which were originally conjectured by Guo. These results establish new $q$-analogues of a supercongruence of Sun.


\vskip 2mm
\noindent{\it MR Subject Classifications}: 33D15, 11A07, 11B65

\section{Introduction}
In recent years, supercongruences on sums involving central binomial coefficients ${2n\choose n}$ have been widely studied by many authors. In 2010, Sun and Tauraso \cite{st-aam-2010} proved that for any odd prime $p$ and positive integer $r$,
\begin{align}
\sum_{k=0}^{p^r-1}\frac{1}{2^k}{2k\choose k}\equiv (-1)^{\frac{p^r-1}{2}}\pmod{p},\label{anew-1}
\end{align}
which was subsequently generalized by Sun \cite{sun-scm-2010} as follows:
\begin{align}
\sum_{k=0}^{p^r-1}\frac{1}{2^k}{2k\choose k}\equiv (-1)^{\frac{p^r-1}{2}}\pmod{p^2}.\label{anew-2}
\end{align}

It is natural to consider $q$-analogues of \eqref{anew-1} and \eqref{anew-2}. To continue the $q$-story of \eqref{anew-1} and \eqref{anew-2}, we recall some $q$-series notation.
The $q$-binomial coefficients are defined as
\begin{align*}
{n\brack k}={n\brack k}_q
=\begin{cases}
\displaystyle\frac{(q;q)_n}{(q;q)_k(q;q)_{n-k}} &\text{if $0\leqslant k\leqslant n$},\\[10pt]
0 &\text{otherwise,}
\end{cases}
\end{align*}
where the $q$-shifted factorials are given by $(a;q)_n=(1-a)(1-aq)\cdots(1-aq^{n-1})$ for $n\ge 1$ and $(a;q)_0=1$. The $q$-integers are defined by $[n]=[n]_q=(1-q^n)/(1-q)$. Moreover, the $n$th cyclotomic polynomial is given by
\begin{align*}
\Phi_n(q)=\prod_{\substack{1\le k \le n\\[3pt](n,k)=1}}
(q-\zeta^k),
\end{align*}
where $\zeta$ denotes a primitive $n$th root of unity.

It is worth mentioning that Guo and Zeng \cite[Corollary 4.2]{gz-aam-2010} first gave an interesting $q$-analogue of \eqref{anew-1} as follows:
\begin{align}
\sum_{k=0}^{n-1}\frac{(q;q^2)_k}{(q;q)_k}q^k\equiv (-1)^{\frac{n-1}{2}}q^{\frac{n^2-1}{4}}\pmod{\Phi_{n}(q)},\label{anew-3}
\end{align}
where $n$ is an odd positive integer. In 2013, Tauraso \cite{tauraso-cm-2013} conjectured that
\eqref{anew-3} also holds modulo $\Phi_{n}(q)^2$, namely,
\begin{align}
\sum_{k=0}^{n-1}\frac{(q;q^2)_k}{(q;q)_k}q^k\equiv (-1)^{\frac{n-1}{2}}q^{\frac{n^2-1}{4}}\pmod{\Phi_{n}(q)^2},\label{anew-4}
\end{align}
which was proved by Guo \cite{guo-ijnt-2019}.
It is clear that letting $n=p^r$ and $q\to 1$ in \eqref{anew-3} and \eqref{anew-4} reduces to
\eqref{anew-1} and \eqref{anew-2}, respectively.

We remark that Wang and Yu \cite{wy-pmh-2022} extended Guo and Zeng's $q$-congruence \eqref{anew-3} as follows:
\begin{align*}
\sum_{k=0}^{n-1}\frac{(q^{2d+1};q^2)_k}{(q;q)_k}q^k\equiv (-1)^{\frac{n-1}{2}+d}q^{\frac{n^2-(2d+1)^2}{4}} \pmod{\Phi_{n}(q)},
\end{align*}
where $n$ is a positive odd integer and $d$ is an integer with $n>2|d|-1$.

Recently, Guo \cite{guo-pmd-2022} established two new $q$-analogues of \eqref{anew-1} as follows:
\begin{align}
\sum_{k=0}^{n-1}\frac{(q;q^2)_k (-1;q^2)_k}{(q^2;q^2)_k}q^{2k}&\equiv (-1)^{\frac{n-1}{2}} \pmod{\Phi_{n}(q)},\label{anew-5}\\[7pt]
\sum_{k=0}^{n-1}\frac{(q;q^2)_k (-q^2;q^2)_k}{(q^2;q^2)_k}q^{2k+1}&\equiv (-1)^{\frac{n-1}{2}} \pmod{\Phi_{n}(q)},\label{anew-6}
\end{align}
where $n$ is a positive odd integer. We remark that the important ingredients in Guo's proof of \eqref{anew-5} and \eqref{anew-6} are a $q$-analogue of Morley's congruence due to Liu, Pan and Zhang \cite{lpz-ade-2015} and Carlitz's identity \cite{carlitz-1974}:
\begin{align}
\sum_{k=0}^n\frac{(a;q)_k (b;q)_k}{(q;q)_k}(-ab)^{n-k}q^{\frac{(n-k)(n+k-1)}{2}}
=\sum_{k=0}^n\frac{(a;q)_{n+1}(-b)^k q^{k\choose 2}}{(q;q)_k (q;q)_{n-k}(1-aq^{n-k})}.
\label{a-3}
\end{align}

Guo \cite[Conjectures 1 and 2]{guo-pmd-2022} also conjectured extensions of \eqref{anew-5} and \eqref{anew-6} as follows:
\begin{conj}
Let $n$ be a positive odd integer. Then, modulo $\Phi_n(q)^2$,
\begin{align}
&\sum_{k=0}^{n-1}\frac{(q;q^2)_k (-1;q^2)_k}{(q^2;q^2)_k}q^{2k}\equiv
\begin{cases}
q^{{n\choose 2}}\quad&\text{if $n\equiv 1\pmod{4}$},\\[7pt]
-q^{{n+1\choose 2}}\quad&\text{if $n\equiv 3\pmod{4}$},
\end{cases}\label{a-1}\\[15pt]
&\sum_{k=0}^{n-1}\frac{(q;q^2)_k (-q^2;q^2)_k}{(q^2;q^2)_k}q^{2k+1}\equiv
\begin{cases}
q^{{n+1\choose 2}}\quad&\text{if $n\equiv 1\pmod{4}$},\\[7pt]
-q^{{n\choose 2}}\quad&\text{if $n\equiv 3\pmod{4}$}.
\end{cases}\label{a-2}
\end{align}
\end{conj}

Note that \eqref{a-1} and \eqref{a-2} give new $q$-analogues of \eqref{anew-2}. Although \eqref{a-1} and \eqref{a-2} are very similar to \eqref{anew-4}, ``the method of proving \eqref{anew-4} given
in \cite{guo-ijnt-2019} does not work here" as mentioned by Guo \cite{guo-pmd-2022}.

Motivated by Guo's proof \cite{guo-pmd-2022}, we shall present a proof of \eqref{a-1} and \eqref{a-2}.
\begin{thm}
The $q$-supercongruences \eqref{a-1} and \eqref{a-2} are true.
\end{thm}

Our proof of \eqref{a-1} and \eqref{a-2} also relies on Carlitz's identity \eqref{a-3}
and a $q$-analogue of Morley's congruence due to Liu, Pan and Zhang \cite{lpz-ade-2015}. We shall prove \eqref{a-1} and \eqref{a-2} in Sections 2 and 3, respectively.

\section{Proof of \eqref{a-1}}
Letting $q\to q^{-1}$ on both sides of \eqref{a-1}, we find that \eqref{a-1} is equivalent to
\begin{align}
\sum_{k=0}^{n-1}\frac{(q;q^2)_k (-1;q^2)_k}{(q^2;q^2)_k}q^{-k^2}\equiv
\begin{cases}
q^{-{n\choose 2}}\quad&\text{if $n\equiv 1\pmod{4}$}\\[7pt]
-q^{-{n+1\choose 2}}\quad&\text{if $n\equiv 3\pmod{4}$}
\end{cases}\pmod{\Phi_n(q)^2}.\label{b-1}
\end{align}

Taking $q\to q^2,a=q,b=-1$ and $n\to n-1$ in \eqref{a-3}, we are led to the following identity (see also \cite[(2.3)]{guo-pmd-2022}):
\begin{align}
&\sum_{k=0}^{n-1}\frac{(q;q^2)_k(-1;q^2)_k}{(q^2;q^2)_k}q^{(n-1)^2-k^2}\notag\\[7pt]
&=\sum_{k=0}^{n-1}\frac{(q;q^2)_n q^{k^2-k}}{(q^2;q^2)_k (q^2;q^2)_{n-k-1}(1-q^{2n-2k-1})}\notag\\[7pt]
&=\sum_{k=0}^{n-1}\frac{(q;q^2)_n q^{k^2-k}}{(q^2;q^2)_{n-1}(1-q^{2n-2k-1})}{n-1\brack k}_{q^2}.\label{b-2}
\end{align}
Let $a_{n,k}$ denote the summand on the right-hand side of \eqref{b-2}, namely,
\begin{align}
a_{n,k}=\frac{(q;q^2)_n q^{k^2-k}}{(q^2;q^2)_{n-1}(1-q^{2n-2k-1})}{n-1\brack k}_{q^2}.\label{bnew-1}
\end{align}

Note that for $0\le k\le n-1$,
\begin{align*}
{n-1\brack k}_{q^2}&=\frac{(1-q^{2n-2})(1-q^{2n-4})\cdots(1-q^{2n-2k})}{(1-q^2)(1-q^4)\cdots(1-q^{2k})}\\[7pt]
&\equiv \frac{(1-q^{-2})(1-q^{-4})\cdots(1-q^{-2k})}{(1-q^2)(1-q^4)\cdots(1-q^{2k})} \pmod{\Phi_n(q)}\\[7pt]
&=(-1)^kq^{-k^2-k}.
\end{align*}
Moreover, the $q$-shifted factorial $(q;q^2)_n$ contains the factor $1-q^n$ which is divisible by $\Phi_n(q)$ and $(q^2;q^2)_{n-1}(1-q^{2n-2k-1})$ is coprime with $\Phi_n(q)$ except for $k=\frac{n-1}{2}$.
It follows that for $k\not=\frac{n-1}{2}$,
\begin{align}
a_{n,k}&\equiv \frac{(-1)^{k} (1-q)(q;q^2)_{n-1} }{(1-q^{2k+1})(q^2;q^2)_{n-1}} \pmod{\Phi_n(q)^2},
\label{b-3}
\end{align}
where we have used the fact that $q^{sn}\equiv 1\pmod{\Phi_n(q)}$ for any integer $s$.

By \cite[(2.4) and (2.7)]{wn-pmh-2022}, we have
\begin{align}
\frac{(q;q^2)_{n-1} }{(q^2;q^2)_{n-1}}=\frac{(q;q^2)_{n-1} }{(q;q)_{n-1} (-q;q)_{n-1}}
\equiv \frac{-q[n]}{(-q;q)_{n-1}}\equiv -q[n]\pmod{\Phi_n(q)^2}.\label{b-4}
\end{align}
Substituting \eqref{b-4} into the right-hand side of \eqref{b-3} gives
\begin{align}
a_{n,k}&\equiv q(1-q^n) \frac{(-1)^{k+1} }{1-q^{2k+1}} \pmod{\Phi_n(q)^2}
\quad\text{for $k\not=\frac{n-1}{2}$}.\label{b-5}
\end{align}
Let $b_{n,k}$ denote the right-hand side of \eqref{b-5}, namely,
\begin{align}
b_{n,k}=q(1-q^n) \frac{(-1)^{k+1} }{1-q^{2k+1}}.\label{b-6}
\end{align}

It follows from \eqref{b-2}, \eqref{bnew-1}, \eqref{b-5} and \eqref{b-6} that
\begin{align}
&\sum_{k=0}^{n-1}\frac{(q;q^2)_k(-1;q^2)_k}{(q^2;q^2)_k}q^{(n-1)^2-k^2}\notag\\[7pt]
&=\sum_{\substack{0\le k \le n-1\\[3pt]k\not= (n-1)/2}} a_{n,k}+a_{n,(n-1)/2}\notag\\[7pt]
&\equiv \sum_{\substack{0\le k \le n-1\\[3pt]k\not= (n-1)/2}} b_{n,k}+a_{n,(n-1)/2}\pmod{\Phi_n(q)^2}.
\label{b-8}
\end{align}

A special case of a $q$-analogue of Morley's congruence \cite[(1.5)]{lpz-ade-2015} (see also \cite[(2.5)]{guo-pmd-2022}) reads
\begin{align}
{n-1\brack \frac{n-1}{2}}_{q^2}\equiv (-1)^{\frac{n-1}{2}}q^{\frac{1-n^2}{4}}(-q;q)_{n-1}^2\pmod{\Phi_n(q)^2}.\label{b-9}
\end{align}
From \eqref{bnew-1} and \eqref{b-9}, we deduce that
\begin{align*}
a_{n,(n-1)/2}&=\frac{q^{\frac{n^2-4n+3}{4}}(q;q^2)_n }{(1-q^{n})(q^2;q^2)_{n-1}}{n-1\brack \frac{n-1}{2}}_{q^2}\\[7pt]
&\equiv\frac{(-1)^{\frac{n-1}{2}}q^{-n+1}(q;q^2)_n (-q;q)_{n-1}^2 }{(1-q^{n})(q^2;q^2)_{n-1}} \pmod{\Phi_n(q)^2}.
\end{align*}
It is easy to check that
\begin{align*}
\frac{(q;q^2)_n (-q;q)_{n-1}^2 }{(1-q^{n})(q^2;q^2)_{n-1}}
={2n\brack n}\frac{1}{1+q^n}.
\end{align*}
Thus,
\begin{align}
a_{n,(n-1)/2}
\equiv \frac{(-1)^{\frac{n-1}{2}}q^{-n+1}}{1+q^n}{2n\brack n}\pmod{\Phi_n(q)^2}.
\label{b-10}
\end{align}
Using the following $q$-congruence (see \cite[(4.4)]{liu-rj-2023}):
\begin{align*}
{2n\brack n}\equiv 2-n(1-q^n)\pmod{\Phi_n(q)^2},
\end{align*}
we obtain
\begin{align}
a_{n,(n-1)/2}&\equiv \frac{(-1)^{\frac{n-1}{2}}q^{-n+1}}{1+q^n}\left(2-n(1-q^n)\right)\notag\\[7pt]
&\equiv (-1)^{\frac{n-1}{2}}\left(q^{-n+1}+\frac{(1-n)q(1-q^n)}{2}\right)
\pmod{\Phi_n(q)^2}.\label{b-11}
\end{align}

Next, we shall determine the following summation modulo $\Phi_n(q)^2$:
\begin{align}
\sum_{\substack{0\le k \le n-1\\[3pt]k\not= (n-1)/2}} b_{n,k}=-q(1-q^n)\sum_{\substack{0\le k \le n-1\\[3pt]k\not= (n-1)/2}} \frac{(-1)^{k} }{1-q^{2k+1}}.\label{b-12}
\end{align}
Note that
\begin{align}
&\sum_{\substack{0\le k \le n-1\\[3pt]k\not= (n-1)/2}} \frac{(-1)^{k} }{1-q^{2k+1}}\notag\\[7pt]
&=\sum_{k=0}^{\frac{n-3}{2}}\frac{(-1)^{k} }{1-q^{2k+1}}+\sum_{k=\frac{n+1}{2}}^{n-1}\frac{(-1)^{k} }{1-q^{2k+1}}\notag\\[7pt]
&=\sum_{k=0}^{\frac{n-3}{2}}\frac{(-1)^{k} }{1-q^{2k+1}}+\sum_{k=0}^{\frac{n-3}{2}}\frac{(-1)^{n-k-1} }{1-q^{2n-2k-1}}\notag\\[7pt]
&\equiv \sum_{k=0}^{\frac{n-3}{2}}\frac{(-1)^{k} }{1-q^{2k+1}}-\sum_{k=0}^{\frac{n-3}{2}}\frac{(-1)^{k}q^{2k+1} }{1-q^{2k+1}}\pmod{\Phi_n(q)}\notag\\[7pt]
&=\sum_{k=0}^{\frac{n-3}{2}}(-1)^k\notag\\[7pt]
&=\frac{1+(-1)^{\frac{n-3}{2}}}{2}.\label{b-13}
\end{align}
It follows from \eqref{b-12} and \eqref{b-13} that
\begin{align}
\sum_{\substack{0\le k \le n-1\\[3pt]k\not= (n-1)/2}} b_{n,k}\equiv -\frac{1+(-1)^{\frac{n-3}{2}}}{2}q(1-q^n)\pmod{\Phi_n(q)^2}.\label{b-14}
\end{align}

Furthermore, combining \eqref{b-8}, \eqref{b-11} and \eqref{b-14} gives
\begin{align}
&\sum_{k=0}^{n-1}\frac{(q;q^2)_k(-1;q^2)_k}{(q^2;q^2)_k}q^{(n-1)^2-k^2}\notag\\[7pt]
&\equiv -\frac{\left(1+(-1)^{\frac{n-3}{2}}\right)q(1-q^n)}{2}\notag\\[7pt]
&+(-1)^{\frac{n-1}{2}}\left(q^{-n+1}+\frac{(1-n)q(1-q^n)}{2}\right)
\pmod{\Phi_n(q)^2}.\label{b-15}
\end{align}

Next, we shall distinguish two cases to prove \eqref{b-1}.

{\noindent \bf Case 1}~~~$n\equiv 1\pmod{4}$. By \eqref{b-15}, we have
\begin{align}
\sum_{k=0}^{n-1}\frac{(q;q^2)_k(-1;q^2)_k}{(q^2;q^2)_k}q^{(n-1)^2-k^2}
\equiv q^{1-n}+\frac{(1-n)q(1-q^n)}{2}\pmod{\Phi_n(q)^2}.\label{b-16}
\end{align}
Since
\begin{align}
q^{sn}&=1-(1-q^{sn})\notag\\[10pt]
&=1-(1-q^n)(1+q^n+q^{2n}+\cdots+q^{(s-1)n})\notag\\[10pt]
&\equiv 1-s(1-q^n)\pmod{\Phi_n(q)^2},\label{b-17}
\end{align}
we have
\begin{align}
q^{(n-1)^2-{n\choose 2}}&=q^{1-n}q^{\frac{n(n-1)}{2}}\notag\\
&\equiv q^{1-n}\left(1-\frac{(n-1)(1-q^n)}{2}\right)\notag\\
&\equiv q^{1-n}+\frac{(1-n)q(1-q^n)}{2}\pmod{\Phi_n(q)^2}.\label{b-18}
\end{align}
Observe that both of the right-hand sides of \eqref{b-16} and \eqref{b-18} are equal.
Then the proof of the case $n\equiv 1\pmod{4}$ of \eqref{b-1} follows from \eqref{b-16} and \eqref{b-18}.

{\noindent \bf Case 2}~~~$n\equiv 3\pmod{4}$. By \eqref{b-15}, we have
\begin{align}
\sum_{k=0}^{n-1}\frac{(q;q^2)_k(-1;q^2)_k}{(q^2;q^2)_k}q^{(n-1)^2-k^2}
\equiv -q^{1-n}+\frac{(n-3)q(1-q^n)}{2}\pmod{\Phi_n(q)^2}.\label{b-19}
\end{align}
Using \eqref{b-17}, we obtain
\begin{align}
-q^{(n-1)^2-{n+1\choose 2}}&=-q^{1-n}q^{\frac{n(n-3)}{2}}\notag\\
&\equiv -q^{1-n}\left(1-\frac{(n-3)(1-q^n)}{2}\right)\notag\\
&\equiv -q^{1-n}+\frac{(n-3)q(1-q^n)}{2}\pmod{\Phi_n(q)^2}.\label{b-20}
\end{align}
Then the proof of the case $n\equiv 3\pmod{4}$ of \eqref{b-1} follows from \eqref{b-19} and \eqref{b-20}.

\section{Proof of \eqref{a-2}}
Letting $q\to q^{-1}$ on both sides of \eqref{a-2}, we find that \eqref{a-2} is equivalent to
\begin{align}
\sum_{k=0}^{n-1}\frac{(q;q^2)_k (-q^2;q^2)_k}{(q^2;q^2)_k}q^{-(k+1)^2}\equiv
\begin{cases}
q^{-{n+1\choose 2}}\quad&\text{if $n\equiv 1\pmod{4}$}\\[7pt]
-q^{-{n\choose 2}}\quad&\text{if $n\equiv 3\pmod{4}$}
\end{cases}\pmod{\Phi_n(q)^2}.\label{c-1}
\end{align}

Letting $q\to q^2, a=q,b=-q^2$ and $n\to n-1$ in \eqref{a-3} reduces to the following identity (see also \cite[(2.8)]{guo-pmd-2022}):
\begin{align}
&\sum_{k=0}^{n-1}\frac{(q;q^2)_k(-q^2;q^2)_k}{(q^2;q^2)_k}q^{n^2-(k+1)^2}\notag\\[7pt]
&=\sum_{k=0}^{n-1}\frac{(q;q^2)_n q^{k^2+k}}{(q^2;q^2)_k (q^2;q^2)_{n-k-1}(1-q^{2n-2k-1})}\notag\\[7pt]
&=\sum_{k=0}^{n-1}\frac{(q;q^2)_n q^{k^2+k}}{(q^2;q^2)_{n-1}(1-q^{2n-2k-1})}{n-1\brack k}_{q^2}.
\label{c-2}
\end{align}
Let $c_{n,k}$ denote the summand on the right-hand side of \eqref{c-2}, namely,
\begin{align}
c_{n,k}=\frac{(q;q^2)_n q^{k^2+k}}{(q^2;q^2)_{n-1}(1-q^{2n-2k-1})}{n-1\brack k}_{q^2}=q^{2k}a_{n,k}.
\label{c-3}
\end{align}

By \eqref{b-5} \eqref{b-6} and \eqref{c-3}, we have
\begin{align}
c_{n,k}&\equiv q(1-q^n) \frac{(-1)^{k+1}q^{2k} }{1-q^{2k+1}}\pmod{\Phi_n(q)^2}\notag\\[7pt]
&=(-1)^k(1-q^n)+q^{-1}b_{n,k}\quad\text{for $k\not=\frac{n-1}{2}$}.\label{c-4}
\end{align}
It follows from \eqref{b-14} and \eqref{c-4} that
\begin{align}
\sum_{\substack{0\le k \le n-1\\[3pt]k\not= (n-1)/2}} c_{n,k}
&\equiv (1-q^n)\sum_{\substack{0\le k \le n-1\\[3pt]k\not= (n-1)/2}}(-1)^k
+q^{-1}\sum_{\substack{0\le k \le n-1\\[3pt]k\not= (n-1)/2}}b_{n,k}\notag\\[7pt]
&\equiv \left(1-(-1)^{\frac{n-1}{2}}\right)(1-q^n)-\frac{1+(-1)^{\frac{n-3}{2}}}{2}(1-q^n)\pmod{\Phi_n(q)^2}\notag\\[7pt]
&=\frac{1-(-1)^{\frac{n-1}{2}}}{2}(1-q^n).\label{c-5}
\end{align}

Furthermore, by \eqref{b-11} and \eqref{c-3} we have
\begin{align}
c_{n,(n-1)/2}&=q^{n-1}a_{n,(n-1)/2}\notag\\[7pt]
&\equiv (-1)^{\frac{n-1}{2}}\left(1+\frac{(1-n)(1-q^n)}{2}\right)
\pmod{\Phi_n(q)^2}.\label{c-6}
\end{align}

It follows from \eqref{c-2}, \eqref{c-5} and \eqref{c-6} that
\begin{align}
&\sum_{k=0}^{n-1}\frac{(q;q^2)_k(-q^2;q^2)_k}{(q^2;q^2)_k}q^{n^2-(k+1)^2}\notag\\[7pt]
&=\sum_{\substack{0\le k \le n-1\\[3pt]k\not= (n-1)/2}} c_{n,k}+c_{n,(n-1)/2}\notag\\[7pt]
&\equiv \frac{1-(-1)^{\frac{n-1}{2}}}{2}(1-q^n)+(-1)^{\frac{n-1}{2}}\left(1+\frac{(1-n)(1-q^n)}{2}\right) \pmod{\Phi_n(q)^2}.\label{c-7}
\end{align}

We shall distinguish two cases to prove \eqref{c-1}.

{\noindent \bf Case 1}~~~$n\equiv 1\pmod{4}$. By \eqref{c-7}, we have
\begin{align}
\sum_{k=0}^{n-1}\frac{(q;q^2)_k(-q^2;q^2)_k}{(q^2;q^2)_k}q^{n^2-(k+1)^2}
\equiv 1+\frac{(1-n)(1-q^n)}{2} \pmod{\Phi_n(q)^2}.\label{c-8}
\end{align}
On the other hand, by \eqref{b-17} we have
\begin{align}
q^{n^2-{n+1\choose 2}}=q^{\frac{n(n-1)}{2}}\equiv 1+\frac{(1-n)(1-q^n)}{2}\pmod{\Phi_n(q)^2}.
\label{c-9}
\end{align}
Combining \eqref{c-8} and \eqref{c-9}, we complete the proof of the case $n\equiv 1\pmod{4}$ of \eqref{c-1}.

{\noindent \bf Case 2}~~~$n\equiv 3\pmod{4}$. In this case, the $q$-congruence \eqref{c-7} reduces to
\begin{align}
\sum_{k=0}^{n-1}\frac{(q;q^2)_k(-q^2;q^2)_k}{(q^2;q^2)_k}q^{n^2-(k+1)^2}
\equiv -1+\frac{(n+1)(1-q^n)}{2}\pmod{\Phi_n(q)^2}.\label{c-10}
\end{align}
By \eqref{b-17}, we have
\begin{align}
-q^{n^2-{n\choose 2}}
&=-q^{\frac{n(n+1)}{2}}\notag\\
&\equiv -1+\frac{(n+1)(1-q^n)}{2}\pmod{\Phi_n(q)^2}.\label{c-11}
\end{align}
Then the proof of the case $n\equiv 3\pmod{4}$ of \eqref{c-1} follows from \eqref{c-10} and \eqref{c-11}.

\vskip 5mm \noindent{\bf Acknowledgments.}
The first author was supported by the National Natural Science Foundation of China (grant 12171370).


\end{document}